\documentclass[a4paper, reqno]{amsart}
\usepackage{graphicx}
\usepackage{amssymb,amsmath}
\usepackage{amsthm}
\usepackage{paralist}
\usepackage{mathtools}
\usepackage{color}
\usepackage[pdfborder={0 0 0}]{hyperref} 
\newenvironment{bew}{\begin{proof}[\textbf{Proof}]}{\end{proof}}

\newtheorem{satz}{Proposition}[section]
\newtheorem{theo}[satz]{Theorem}
\newtheorem{lemma}[satz]{Lemma}
\theoremstyle{definition}
\newtheorem{defi}[satz]{Definition}
\newtheorem{remark}[satz]{Remark}

\numberwithin{equation}{section}

\renewcommand{\epsilon}{\varepsilon}
\renewcommand{\Im}{\operatorname{Im}}

\newcommand{\tr}{\operatorname{tr}}
\newcommand{\sgn}{\operatorname{sgn}}

\newcommand{\spann}{\operatorname{span}_\R}
\newcommand\ko{\hspace{.1em}}

\newcommand{\R}{\mathbb{R}}
\newcommand{\C}{\mathbb{C}}
\newcommand{\Z}{\mathbb{Z}}
\newcommand{\N}{\mathbb{N}}
\renewcommand{\H}{\mathbb{H}}
\newcommand{\myhomopol}[1]{\mathcal{P}_{#1}^{m,n}}
\newcommand{\comp}{{C_Q}}

\newcommand{\Amaj}{M}
\newcommand{\mymatrix}[2]{{#2}^{{#1}\times {#1}}}
\newcommand{\Aneg}{A^{-}}
\newcommand{\Apos}{A^{+}}
\newcommand{\Usenk}{U^\perp}
\newcommand{\Uc}{U^c}
\newcommand{\trans}{\mathsf{T}\hspace{-2pt}}

\newcommand{\uvector}{\boldsymbol{u}}
\newcommand{\vvector}{\boldsymbol{v}}
\newcommand{\xvector}{\boldsymbol{x}}
\newcommand{\kvector}{\boldsymbol{k}}
\newcommand{\tvector}{\boldsymbol{t}}
\newcommand{\ucolumnvector}[1]{\boldsymbol{u_{#1}}}
\newcommand{\mycomponentvectorgeneral}{\boldsymbol{c}}
\newcommand{\mycomponentvector}[1]{\boldsymbol{c_{#1}}}
\newcommand{\mysamesign}[1]{\widetilde{x}_{#1}}
\newcommand{\mysamesigncolumnvector}[1]{\boldsymbol{x_{#1}}}

\newcommand{\usenk}{\boldsymbol{u^\perp}}
\newcommand{\ucolumnvectorsenk}[1]{\boldsymbol{u_{#1}^\perp}}

\newcommand{\ucolumnvectorc}[1]{\boldsymbol{u_{#1}^c}}

\newcommand{\Quadratic}{\boldsymbol{Q}}
\newcommand{\Quadraticplus}{\boldsymbol{Q}^{+}}
\newcommand{\Quadraticminus}{\boldsymbol{Q}^{-}}

\newcommand{\mytheta}[1]{\vartheta_{#1}}

\newcommand{\laplace}{\boldsymbol{\Delta}}
\newcommand{\laplacian}{\operatorname{tr}\boldsymbol{\Delta}}

\newcommand{\partialmatrix}[1]{\frac{\partial}{\partial #1}}

\newcommand{\Qkl}{Q^{+}_{k,\ell}}

\newenvironment{acknowledgments}{
   \abstract}{
  \endabstract
}
\makeatletter
\@namedef{subjclassname@2020}{\textup{2020} Mathematics Subject Classification}
\makeatother

\title{Siegel theta series for quadratic forms of signature $(m-1,1)$}
\author{Christina Roehrig}
\address{University of Cologne, Department of Mathematics and Computer Science, Division of Mathematics, Weyertal 86-90, 50931 Cologne, Germany}
\email{croehrig@math.uni-koeln.de}
\begin{document}
\begin{abstract}
We investigate Siegel theta series for quadratic forms of signature $(m-1,1)$. On the one hand, we construct a holomorphic series that does not transform like a modular form. On the other hand, we construct a non-holomorphic series that transforms like a Siegel modular form of weight $m/2$. Moreover, the holomorphic series describes almost everywhere the holomorphic part of the modular series.
\end{abstract}
\maketitle
\section{Introduction}
While Siegel modular forms play an important role in various areas of mathematics, such as number theory and algebraic geometry, the number of explicit constructions is rather limited. We can obtain interesting examples by considering Siegel theta series that are associated with quadratic forms $\Quadratic(U)=\frac12\tr (U^\trans A U)$ on $\R^{m\times n}$, where $A\in \Z^{m\times m}$ is an even symmetric and non-degenerate matrix with signature $(r,s)$.

If $A$ is positive definite, we can construct Siegel theta series as follows (see \cite{roehrig} for a more detailed description). We consider functions of the form
$p=\exp ( -\laplacian_{A}/8\pi)P$,
where 
\[\exp \Bigl( -\frac{\laplacian_{A}}{8\pi}\Bigr)=\sum_{k=0}^\infty \frac{(-1)^k}{(8\pi)^k\ko k!}\ko(\laplacian_A)^k\quad\text{with }\laplace_A=\Bigl(\partialmatrix{U}\Bigr)^\trans\ko A^{-1}\partialmatrix{U},\]
and $P:\C^{m\times n}\longrightarrow\C$ is a polynomial with the homogeneity property $P(UN)=\det N^\alpha P(U)$ for all $N\in \C^{n\times n}$ and a fixed $\alpha\in \N_0$.
We define the Siegel theta series associated with $p$ as
\[\mytheta{p}(Z)=\det Y^{-\alpha/2}\sum\limits_{U\in\Z^{m\times n}} p(UY^{1/2})\ko\exp \bigl(\pi i \tr (U^\trans AUZ)\bigr)\quad(\H_n \ni Z=X+iY),\]
where $\H_n$ is the Siegel upper half-space and $Y^{1/2}$ is the square root of the positive definite matrix $Y$. Then $\mytheta{p}$ transforms like a Siegel modular form of weight $m/2+\alpha$.
As the homogeneity property of $P$ is not maintained when we apply the operator $\exp ( -\laplacian_{A}/8\pi)$, these examples are in general non-holomorphic. 
Only when we consider harmonic polynomials, i.\,e. $(\laplacian_{A})P=0$, we obtain holomorphic functions (see for example Freitag's exposition \cite{Fre83}).

For indefinite quadratic forms it is generally difficult to construct Siegel theta series that are holomorphic \emph{and} modular. Kudla \cite{kudla81}, however, considered quadratic forms of signature $(n,1)$ to construct holomorphic Siegel modular forms of genus $n$ and weight $(n+1)/2$ as the integrals of non-holomorphic theta series. We employ a different approach and recall that we can generalize the construction above to obtain Siegel theta series for indefinite quadratic forms that transform like modular forms (see Borcherds' construction \cite{borcherds} for $n=1$ and our own recent result \cite{roehrig} for higher $n$). In general, these are non-holomorphic and choosing a harmonic polynomial $P$ does not suffice to establish holomorphicity.

In the present paper, we will deal with quadratic forms of signature $(m-1,1)$ and find for arbitrary genus $n\in \N$ holomorphic Siegel theta series that are related to non-holomorphic modular Siegel theta series. We obtain this construction by generalizing the result on elliptic modular forms by Zwegers \cite{zwegers}.
So we review some results on elliptic theta series first.

For positive definite quadratic forms, we have Schoeneberg's description \cite{schoeneberg} for the case that $m$ is even and the result by Shimura \cite{shimura} for the case that $m$ is odd. For $Q(\uvector)=\frac12 \uvector^\trans A \uvector$, where $A$ has signature $(m-1,1)$, G\"ottsche and Zagier \cite{GZ98} introduced holomorphic theta series, which were then modified by Zwegers \cite{zwegers} to construct theta series with modular transformation behavior. Instead of defining the theta series as a series over a full lattice in $\R^m$, one sums over a suitable cone in $\R^m$ to ensure the absolute convergence of the series. We introduce some notation to make this more explicit. We refer to one of the components of $\{\mycomponentvectorgeneral\in \R^m \mid Q(\mycomponentvectorgeneral)<0\}$ as $C_Q$ and define the theta function depending on two vectors $\mycomponentvector{0},\mycomponentvector{1}\in C_Q$. Further, let $\boldsymbol{h},\kvector\in \R^m$, and let $B$ denote the bilinear form associated with $Q$.
Then the holomorphic theta series is given by
\[\vartheta_{\mycomponentvector{0},\mycomponentvector{1}} (z)= \sum_{\uvector\in \boldsymbol{h}+\Z^m}\bigl\{\sgn \bigl(B(\mycomponentvector{0},\uvector)\bigr) - \sgn \bigl(B(\mycomponentvector{1},\uvector)\bigr)\bigr\}\ko \exp\bigl(2\pi i Q(\uvector)z+2\pi i B(\uvector,\kvector)\bigr)\quad(z \in \H)\]
and the modular theta series by
\[\widehat\vartheta_{\mycomponentvector{0},\mycomponentvector{1}} (z)= \sum_{\uvector\in\boldsymbol{h}+ \Z^m}\Bigl\{E\Bigl(\frac{B(\mycomponentvector{0},\uvector)}{\sqrt{-Q(\mycomponentvector{0})}}y^{1/2}\Bigr)-E\Bigl(\frac{B(\mycomponentvector{1},\uvector)}{\sqrt{-Q(\mycomponentvector{1})}}y^{1/2}\Bigr)\Bigr\}\ko \exp \bigl(2\pi i Q(\uvector)z+2\pi i B(\uvector,\kvector)\bigr),\]
where $y=\Im z>0$ and
\[E(x) = 2\int_0^x \exp(-\pi v^2)\ko dv = \sgn (x) - \sgn (x) \int_{x^2}^\infty v^{-1/2} \exp(-\pi v)\ko dv\quad (x\in \R).\]
Since $E(xy^{1/2})\rightarrow \sgn (x)$ for $y\rightarrow \infty$, the theta series $\vartheta_{\mycomponentvector{0},\mycomponentvector{1}}$ describes the holomorphic part of  $\widehat\vartheta_{\mycomponentvector{0},\mycomponentvector{1}}$ and therefore we call $\widehat\vartheta_{\mycomponentvector{0},\mycomponentvector{1}}$ the modular completion of $\vartheta_{\mycomponentvector{0},\mycomponentvector{1}}$.

For signature $(m-2,2)$, holomorphic theta series and their modular completions were constructed in the work of Alexandrov, Banerjee, Manschot, and Pioline \cite{ABMP18}. Their suggestion for a generalization to generic signature $(r,s)$ was then explicitly carried out by Nazaroglu \cite{Naz18}. In a slightly different setting, considering positive polyhedral cones, similar results were presented by Raum \cite{raum}. Also Funke and Kudla \cite{FK17,FK19} gave a general framework for the construction of non-holomorphic theta series for indefinite quadratic forms.

From a geometric point of view, an analogue of these theta series can be constructed as integrals of the theta forms introduced by Kudla and Millson \cite{KM86,KM87,KM90} (in the latter, Siegel modular forms of higher genus are also discussed). The connection between the geometric and the classical approach was established in the aforementioned works by Funke and Kudla \cite{FK17,FK19} and in the special cases where $s=1,2$ by Kudla \cite{kudla13,kudla18}.

We will employ this connection in the following 
when we describe modular Siegel theta series for quadratic forms of signature $(m-1,1)$. This has also been done by Livinsky \cite{livinsky}. He deduces the connection to Zwegers' theta series, gives a construction for arbitrary $n\in \N$ and furthermore a very explicit description of these theta series for the case $n=2$. We obtain almost the same construction for the modular version of the theta series.  However, a connection to holomorphic series is not given in \cite{livinsky}. In contrast to that, we determine an explicit construction of the associated holomorphic Siegel theta series and establish a connection between the holomorphic and the modular version. Note that we have to assume $m>n$ to obtain non-vanishing series.

We state the main results in more detail in the next section. There we also introduce the notation and several important definitions that will be used throughout the rest of the paper. Besides, we shortly recapitulate a result of \cite{roehrig}, which states that theta series arising from certain functions transform like Siegel modular forms. In Section \ref{section_holomorphic}, we construct holomorphic Siegel theta series and thus prove part (i) of the main theorem. In Section \ref{section_modular}, we consider Siegel theta series with modular transformation behavior (this is part (ii) of the main theorem). To conclude this section, we show that the function used in the construction of the modular series asymptotes to the function employed in the construction of the holomorphic one, which proves part (iii). We will also see there that for $n=1$ we get back Zwegers' result \cite{zwegers}.
\section{Definitions, previous results and statement of the main results}
We define the Siegel upper half-space as 
\[
\H_n :=\lbrace Z=X+iY \mid X,Y \in \mymatrix{n}{\R}\text{ symmetric and } Y \text{ positive definite}\rbrace
\]
and the full Siegel modular group $$\Gamma_n:=\big\lbrace M=\left(\begin{smallmatrix}
A&B\\
C&D
\end{smallmatrix}\right)\in \mymatrix{2n}{\Z} \mid M^\trans JM=J \big\rbrace,\quad\text{where}\quad J=\left(\begin{smallmatrix}
\mathrm{O}&I_n\\
-I_n&\mathrm{O}
\end{smallmatrix}\right),$$ 
which operates on $\H_n$ by
$$Z\mapsto M\langle Z\rangle=(AZ+B)(CZ+D)^{-1}.$$
The group $\Gamma_n$ is generated by the matrices $
\left(\begin{smallmatrix}
I_n&S\\
\mathrm{O}&I_n
\end{smallmatrix}\right)$, where $S\in \Z^{n\times n}$ is symmetric, and the matrix
$\left(\begin{smallmatrix}
\mathrm{O}&-I_n\\
I_n&\mathrm{O}
\end{smallmatrix}\right),$ cf. \cite[p.\ko322-328]{Fre83}.
We define Siegel theta series of the following form:
\begin{defi}
Let $\H_n \ni Z=X+iY$ and let $A\in \mymatrix{m}{\Z}$ be an even symmetric and non-degenerate matrix with signature $(r,s)$. The theta series with characteristics $\mathcal{H,K}\in \R^{m\times n}$, associated with $p:\R^{m\times n}\longrightarrow\R$ and $A$ is defined as
\begin{align*}
\mytheta{p}(Z)=\vartheta_{\mathcal{H},\mathcal{K},p,A}(Z):=\sum\limits_{U\in \mathcal{H}+ \Z^{m\times n}} p(UY^{1/2})\ko\exp \bigl(\pi i \tr (U^\trans AUZ)+2\pi i\tr(\mathcal{K}^\trans AU)\bigr).
\end{align*}
\end{defi}
\begin{remark}
(a) We have to ensure by the choice of $p$ that the theta series $\mytheta{p}$ is defined by an absolutely convergent series (which is not obvious as $A$ is indefinite). For the modular as well as for the holomorphic version, we will see that this is satisfied (see Remark \ref{rem_convergence} and Proposition \ref{proposition_holomorphic}).\\
(b) Note that we normally also have the factor $\det Y^{-\lambda/2}$ in the definition of $\mytheta{p}$, where the choice of $\lambda\in \Z$ depends on $p$, but here we will consider theta series with $\lambda=0$.\\
(c) As we take $A$ to be fixed, we usually drop this parameter in the index. We do the same for the characteristics $\mathcal{H,K}$, as they only play a role when we determine the explicit modular transformation behavior. 
\end{remark}

In \cite{roehrig} we specifically constructed modular Siegel theta series of this form. In order to do so, we split up $A$ into a positive semi-definite part $\Apos$ and a negative semi-definite part $\Aneg$. Further, we define $M$ as a positive definite majorant matrix associated with $A$, i.\,e. $M=\Apos -\Aneg$.
Then we can also write $U\in \R^{m\times n}$ as $U=U^++U^-$, where $U^+$ lies in a subspace of $\R^{m\times n}$ on which the quadratic form is positive semi-definite and $U^-$ in a subspace of $\R^{m\times n}$ on which the quadratic form is negative semi-definite. 
In \cite{roehrig} we obtained the decomposition of $A$ by considering the eigenvectors of $A$. However, we can consider any decomposition of this form, see \eqref{align_split} for the one used in the present paper.

For $\alpha\in \N_0$, let $\myhomopol{\alpha}$ denote the vector space of polynomials $P:\R^{m\times n}\longrightarrow\R$ with the homogeneity property $P(UN)=\det N^\alpha P(U)$ for all $N\in\R^{n\times n}$. Further, for $M\in \Z^{\mu\times \mu}$ we define $M_0\in \Z^{\mu\times \mu}$ by $(M_0)_{ij}=M_{ii}$ for $i=j$ and zero otherwise and we use the notation $\mathrm{1}_{nm}$ for a matrix with $n$ rows and $m$ columns, whose entries are all equal to 1.

A special case of Lemma 4.2 and Proposition 4.10 in \cite{roehrig} is the following result, where we choose $\alpha,\beta\in \N_0$ with $\alpha-\beta=s$ so that the weight is $m/2$. Also note that in contrast to \cite{roehrig} we choose $A$ to be even in the present paper.
\begin{theo}\label{theorem_solutions}
Let $\alpha,\beta\in \N_0$ with $\alpha-\beta=s$. Further, let $P$ be defined as the product $P(U)=P_r(U^+)\cdot P_s(U^-)$ with $P_r\in \myhomopol{\alpha}$ and $P_s\in \myhomopol{\beta}$ and set
\begin{align}\label{align_coefficient_indefinite}
p(U)=\exp\Bigl( -\frac{\laplacian_{\Amaj}}{8\pi}\Bigr)\bigl(P(U)\bigr)\ko\exp\bigl(2\pi \tr(U^\trans \Aneg U )\bigr).
\end{align}
Then the transformation behavior of $\mytheta{p}$ is as follows. For any symmetric matrix $S\in \Z^{n\times n}$ we have
\begin{align*}
\mytheta{\mathcal{H},\mathcal{K},p,A}(Z+S)=\exp\bigl(-\pi i \tr (\mathcal{H}^\trans A\mathcal{H}S)-\pi i \tr (S_0\mathrm{1}_{nm}A_0\mathcal{H})\bigr)\ko\mytheta{\mathcal{H},\widetilde{\mathcal{K}},p,A}(Z)\quad\text{with }\widetilde{\mathcal{K}}:=\mathcal{K}+\mathcal{H}S,
\end{align*}
and we have
\begin{multline*}
\mytheta{\mathcal{H},\mathcal{K},p,A}(-Z^{-1})=i^{-mn/2}(-1)^{(s/2+\beta)n+\beta s}|\det A|^{-n/2}\det Z^{m/2}\exp \bigl(2\pi i \tr (\mathcal{H}^\trans A \mathcal{K})\bigr)\\
\cdot\sum\limits_{\mathcal{J}\in A^{-1}\Z^{m\times n} \operatorname{mod} \Z^{m\times n} } \mytheta{\mathcal{J}+\mathcal{K},-\mathcal{H},p,A}(Z).
\end{multline*}
\end{theo}
\begin{remark}
We can then either take rational matrices $\mathcal{H}$ and $\mathcal{K}$ and consider these theta series as entries of vector-valued Siegel modular forms or we set $\mathcal{H}=\mathcal{K}=\mathrm{O}$ and obtain scalar-valued modular forms for a certain character and on a suitable congruence subgroup of level $N$ in $\Gamma_n$ (where $N$ is the level of $A$, i.\,e. the smallest $N\in \N$ such that $NA^{-1}$ is an even matrix). However, we restrict ourselves to the description for the generating matrices of $\Gamma_n$.
To put it short, we just say that a Siegel theta series of this kind \emph{transforms like a (Siegel) modular form}.
\end{remark}
\begin{remark}\label{rem_convergence}
For $p:\R^{m\times n}\longrightarrow\R$ as in (\ref{align_coefficient_indefinite}), the series that defines $\mytheta{p}$
is absolutely convergent, because the part of the expression that determines the growth is
\begin{align*}
\exp\bigl(2\pi \tr(U^\trans \Aneg U Y)-\pi \tr(U^\trans AUY)\bigr)=\exp\bigl(-\pi \tr(U^\trans MU Y)\bigr),
\end{align*}
where $M=\Apos-\Aneg$ is the positive definite majorant associated with $A$ described above and $Y$ is positive definite.
\end{remark}
From now on we take the signature of $A$ to be  $(m-1,1)$. This and more notational conventions are fixed in the following definition.
\begin{defi}
Let $A$ be an even symmetric and non-degenerate matrix of signature $(m-1,1)$. Then we define the quadratic forms $\Quadratic:\R^{m\times n}\longrightarrow\R,\ko\Quadratic (U):=\frac12\tr (U^\trans A U)$ and $Q:\R^m\longrightarrow \R,\ko Q(\boldsymbol{u}):=\frac12 \boldsymbol{u}^\trans A\boldsymbol{u}$ with the associated bilinear form $B(\uvector,\vvector)=Q(\uvector+\vvector)-Q(\uvector)-Q(\vvector)$.
\end{defi}
We write henceforth $\ucolumnvector{j}$ for the $j$-th column vector of $U=(\ucolumnvector{1}\ldots \ucolumnvector{n})\in \R^{m\times n}$. We can thus write $\Quadratic$ as $\Quadratic (U)=\sum_{j=1}^n Q(\ucolumnvector{j})$ to consider the quadratic forms on the column vectors of $U$.

Further, we fix an element $\mycomponentvectorgeneral\in \R^m$ with $Q(\mycomponentvectorgeneral)<0$ to split $A$ into a negative semi-definite and a positive semi-definite part. In order to do so, we set
\begin{align}\label{align_split}
\Aneg:=\frac{A\mycomponentvectorgeneral\mycomponentvectorgeneral^\trans A}{2Q(\mycomponentvectorgeneral)}\quad\text{and}\quad \Apos:= A-\Aneg
\end{align}  and we define the corresponding quadratic forms
$\Quadraticminus(U):=\frac12\tr(U^\trans \Aneg U)$ and $\Quadraticplus(U):=\frac12\tr(U^\trans \Apos U)$. In Lemma \ref{lemma_split}, we will see that indeed $\Quadraticminus$ is a negative semi-definite and $\Quadraticplus$ a positive semi-definite quadratic form. We will also show there that if we write $U=\Usenk + \Uc$ by setting $\ucolumnvectorc{j}:=\frac{B(\mycomponentvectorgeneral,\ucolumnvector{j})}{2Q(\mycomponentvectorgeneral)}\mycomponentvectorgeneral$ and $\ucolumnvectorsenk{j}:=\ucolumnvector{j}-\ucolumnvectorc{j}$, the part $\Usenk$ lies in the subspace where $\Quadratic$ is positive semi-definite and $\Uc$ in the subspace where $\Quadratic$ is negative semi-definite.

We recall the construction of elliptic theta series in \cite{zwegers} for quadratic forms of signature $(m-1,1)$. One fixes one of the components in $\R^m$, where $Q$ is negative, by taking a vector $\mycomponentvector{0}\in \R^m$ with $Q(\mycomponentvector{0})<0$ and setting 
\[\comp :=\lbrace \uvector\in \R ^m\mid Q(\uvector)<0,\ko B(\uvector,\mycomponentvector{0})<0\rbrace.\]
The theta series in \cite{zwegers} then depend on two vectors $\mycomponentvector{0},\mycomponentvector{1}\in \comp$.
To generalize this construction to higher genus $n$ we need to define several similar objects:
\begin{defi}
We fix $n+1$ vectors in $\R^m$, which lie in $\comp$, and collect them in a matrix
\begin{align*}
C:=(\mycomponentvector{0}\,\mycomponentvector{1} \ldots \mycomponentvector{n})\in \R^{m\times (n+1)}\quad\text{with }\mycomponentvector{i}\in \comp \subset \R^m.
\end{align*}
We will also consider the matrices
\begin{align*}
\widetilde{C}_i:=(\mycomponentvector{0}\ldots \widehat{\mycomponentvector{i}} \ldots \mycomponentvector{n})\in \R^{m\times n}\quad(0\leq i\leq n),
\end{align*}
where $\widehat{\cdot}$ means that the respective column is omitted. For $U\in \R^{m\times n}$, let $\R^n\ni\mysamesigncolumnvector{i}:=U^\trans A \mycomponentvector{i}$ for $0\leq i \leq n$ and define $\mysamesign{i}$ by setting
\begin{align}\label{align_x_k}
\mysamesign{i}:=(-1)^{i}\det \bigl(\mysamesigncolumnvector{0}\ldots\widehat{\mysamesigncolumnvector{i}}\ldots \mysamesigncolumnvector{n}\bigr)=(-1)^{i}\det \bigl(U^\trans A \widetilde{C}_i\bigr)\quad (0\leq i \leq n).
\end{align}
\end{defi}
We use these $\mysamesign{i}$ to define an absolutely convergent and holomorphic theta series. For the definition of a modular version we consider an $n$-simplex in $\R^m$ that is defined as follows:
\begin{defi}
We define the $n$-simplex
\begin{align*}
S_n:=\bigg\lbrace \sum_{i=0}^n t_i\mycomponentvector{i}\hspace{3pt}\Big\vert\hspace{3pt} t_i\geq 0\text{ for all }0\leq i \leq n\quad\text{and}\quad \sum_{i=0}^n t_i=1\bigg\rbrace.
\end{align*}
\end{defi}
\begin{remark}\label{rem_Sn}
(a) Up to a sign the integral of a certain $n$-form over $S_n$ is independent of the explicit parameterization of $S_n$. So we will only fix a parameterization when we explicitly evaluate the integral.\\
(b) We have $S_n\subset C_Q$: By definition we can write any $\mycomponentvectorgeneral\in S_n$ as $\mycomponentvectorgeneral=\sum_{i=0}^n t_i\mycomponentvector{i}$ with $t_i\geq 0$ for all $0\leq i \leq n$ and $\sum_{i=0}^n t_i=1$. As not all $t_i$ vanish, we have
\[Q(\mycomponentvectorgeneral)=\frac12\sum_{i,j=0}^n t_it_j\ko B(\mycomponentvector{i},\mycomponentvector{j})<0\quad\text{and}\quad B(\mycomponentvectorgeneral,\mycomponentvector{0})=\sum_{i=0}^n t_i\ko B(\mycomponentvector{i},\mycomponentvector{0})<0,\] thus $\mycomponentvectorgeneral\in C_Q$.\\
(c) $S_n$ is a compact set in $\R^m$, because it is closed and bounded.
\end{remark}
We establish a connection between the holomorphic and the modular versions of the theta series. For this purpose we need to transfer the concept of a modular completion that we know for elliptic modular forms to higher genus $n$.
\begin{defi}\label{defi_infinity}
Let $Y\in \R^{n\times n}$ denote a positive definite symmetric matrix. Then all diagonal entries of $Y$ are positive. When all the entries on the diagonal simultaneously go to infinity, we define this as $Y\rightarrow\infty$. If we have a modular theta series $\mytheta{g}$ and a holomorphic theta series $\mytheta{f}$ for which $g(UY^{1/2})\rightarrow f(U)$ for $Y\rightarrow \infty$ holds, we then say that $\mytheta{f}$ describes the holomorphic part of $\mytheta{g}$ and on the other hand $\mytheta{g}$ is referred to as the modular completion of $\mytheta{f}$.
\end{defi}
The main theorem we are going to prove is the following:
\begin{theo}\label{maintheorem}
(i) For
\[
f(U)=f^C(U):=\prod\limits_{i=0}^n\frac{1+\sgn (\mysamesign{i})}{2}-\prod\limits_{i=0}^n\frac{1-\sgn (\mysamesign{i})}{2},
\]
the theta series $\mytheta{f}$ is absolutely convergent and holomorphic in $Z\in \H_n$.\\
(ii) Let
\begin{align*}
g(U)=g^C(U):=\int\limits_{S_n}(-Q(\mycomponentvectorgeneral))^{-n/2}\ko\exp\bigl(2\pi \tr (U^\trans \Aneg U)\bigr)\ko\bigwedge_{j=1}^n  B(\ucolumnvectorsenk{j},d\mycomponentvectorgeneral).
\end{align*}
Then the theta series $\mytheta{g}$ transforms like a Siegel modular form of weight $m/2$.\\
(iii) We have $g(UY^{1/2})\rightarrow f(U)$ almost everywhere for $Y\rightarrow \infty$. 
\end{theo}
\begin{remark}
(a) We define the function $f$ in order to describe the holomorphic part of $g$ as good as possible with a relatively simple function. By definition, $f$ is locally constant and evaluates to $\pm 1$ or $0$ almost everywhere. We see in Proposition \ref{proposition_completion} that the same holds for $g(UY^{1/2})$ for $Y\rightarrow\infty$.
Further, $f(-U)=(-1)^n f(U)$, as $\sgn\bigl(\det ((-U)^\trans A \widetilde{C}_i)\bigr)=(-1)^n \sgn\bigl(\det (U^\trans A \widetilde{C}_i)\bigr)$ for $U\in \R^{m\times n}$, and we also have $g(-U)=(-1)^n g(U)$, which can be deduced immediately from the definition of $g$. Although $f$ and $g$ occur to have similar properties, $\mytheta{f}$ does not exactly describe the holomorphic part of $\mytheta{g}(Z)$ for $Y\rightarrow\infty$.
In Remark \ref{rem_solidangle}, we outline a more precise description of the holomorphic part of $\mytheta{g}$.\\
(b) In order to obtain non-vanishing functions $f$ and $g$, we must assume $m>n$, so we make this assumption throughout the rest of this paper. From Lemma \ref{lemma_fullrank}, it immediately follows that this is a necessary condition so that  $f$ does not vanish identically. The same holds for $g$: the column vectors $\ucolumnvectorsenk{j}$ of $\Usenk$ lie in an $(m-1)$ -  dimensional  subspace of $\R^m$, so for $n\geq m$ we have for all $U\in \R^{m\times n}$ linear dependencies among the $n$ vectors $\ucolumnvectorsenk{j}$. Using the fact that the wedge product is a distributive and alternating map, we deduce that the $n$-form $\bigwedge_{j=1}^n  B(\ucolumnvectorsenk{j},d\mycomponentvectorgeneral)$ vanishes identically and thus in particular $g$.
\end{remark}
\section{Holomorphic Siegel theta series}\label{section_holomorphic}
In this section, we construct theta series of genus $n$ associated with indefinite quadratic forms of signature $(m-1,1)$ that are holomorphic. 
For this purpose, we consider the locally constant functions $f$ described in Theorem \ref{maintheorem}(i).
As in Zwegers' work \cite{zwegers}, we show that the choice of $f$ restricts the summation in the definition of the theta function $\mytheta{f}$ to a component in $\R^{m\times n}$ on which the indefinite form is bounded from below by a positive definite quadratic form.

For $\mysamesign{i}$ as in \eqref{align_x_k}, we define this component as
\begin{multline*}
\mathcal{C}_A:=\lbrace U\in \R^{m\times n} \mid \mysamesign{i}\geq 0\text{ for all }0\leq i \leq n\quad \text{or}\quad \mysamesign{i}\leq 0\text{ for all }0\leq i \leq n,\\
\text{where in both cases not all }\mysamesign{i}\text{ vanish} \rbrace.
\end{multline*}

By definition, $f(U)=0$ for $U\notin \mathcal{C}_A$, as either $\mysamesign{i}=0$ for all $0\leq i \leq n$, or both summands vanish since there exist $i,j\in \lbrace 0,\ldots, n\rbrace$ with $\mysamesign{i}>0$ and $\mysamesign{j}<0$, so the support of $f$ lies in $\mathcal{C}_A$.
Even if $\mathcal{C}_A$ is large, the corresponding theta series might vanish but at least we can exclude choices of $C$ for which $\mathcal{C}_A$ is the empty set and thus $f\equiv 0$. For $n=1$, this is done by choosing two linearly independent vectors $\mycomponentvector{0}$ and $\mycomponentvector{1}$. For higher genus $n$, we show in the next lemma that the matrix $C$ should have full rank, i.\,e. the $n+1$ column vectors are linearly independent (as we assume $m>n$, we can always choose that many linearly independent vectors).
\begin{lemma}\label{lemma_fullrank}
If $C=(\mycomponentvector{0}\,\mycomponentvector{1} \ldots \mycomponentvector{n})$ does not have full rank, $\mathcal{C}_A$ is empty.
\end{lemma}
\begin{proof}
If $\mycomponentvector{0},\mycomponentvector{1}, \ldots, \mycomponentvector{n}\in \R^m$ are linearly dependent, we can write without loss of generality
$\mycomponentvector{0}=\sum_{i=1}^n\lambda_i\mycomponentvector{i}$ for $\lambda_i\in \R$.
We determine $\mysamesign{k}$ for $k\in \{1,\ldots,n\}$, substituting $\mysamesigncolumnvector{0}=\sum_{i=1}^n\lambda_i\mysamesigncolumnvector{i}$:
\[
\mysamesign{k}=(-1)^k\ko\sum_{i=1}^n\lambda_i \det\bigl(\mysamesigncolumnvector{i}\,\mysamesigncolumnvector{1}\ldots\widehat{\mysamesigncolumnvector{k}}\ldots \mysamesigncolumnvector{n}\bigr)
=(-1)^k \lambda_k\ko \det\bigl(\mysamesigncolumnvector{k}\,\mysamesigncolumnvector{1}\ldots\widehat{\mysamesigncolumnvector{k}}\ldots \mysamesigncolumnvector{n}\bigr)
=-\lambda_k\mysamesign{0}.
\]
For $U\in \mathcal{C}_A$, $\mysamesign{0}$ and $\mysamesign{k}$ are both non-positive or both non-negative, so we either have $\lambda_k\leq 0$ or $\mysamesign{0}=\mysamesign{k}=0$. If the latter case holds for any $k$, all $\mysamesign{i}\,(0\leq i\leq n)$ vanish, which contradicts our definition of $\mathcal{C}_A$.
So we have $\lambda_k\leq 0$ for all $k\in \{1,\ldots,n\}$. But since \[0>B(\mycomponentvector{0},\mycomponentvector{0})=\sum_{i=1}^n\lambda_i B(\mycomponentvector{0},\mycomponentvector{i}),\quad\text{where }B(\mycomponentvector{0},\mycomponentvector{i})<0\text{ for all }1\leq i\leq n,\]
not all $\lambda_k$ can be non-positive.
Hence, $\mathcal{C}_A$ is the empty set.
\end{proof}
 
We use the following statement to construct positive definite quadratic forms:
\begin{lemma}[{\cite[Lemma 2.6]{zwegers}}]\label{lemma_pos_def}
Let $\mycomponentvector{k},\mycomponentvector{\ell} \in \comp$ be linearly independent. The quadratic form \[\Qkl:\mathbb{R}^m\longrightarrow \mathbb{R},\quad\Qkl(\boldsymbol{v}):=Q(\boldsymbol{v})+\frac{B(\mycomponentvector{k},\mycomponentvector{\ell})B(\mycomponentvector{k},\boldsymbol{v})B(\mycomponentvector{\ell},\boldsymbol{v})}{4Q(\mycomponentvector{k})Q(\mycomponentvector{\ell})-B(\mycomponentvector{k},\mycomponentvector{\ell})^2}\] is positive definite.
\end{lemma}
For $U\in\mathcal{C}_A$ we have:
\begin{lemma}\label{lemma_sign}
For any column vector $\ucolumnvector{j}$ of $U\in \mathcal{C}_A$, there exist $k,\ell\in\lbrace 0,\ldots,n\rbrace$ such that \[\sgn\bigl(B(\mycomponentvector{k},\ucolumnvector{j})\bigr)\neq \sgn\bigl(B(\mycomponentvector{\ell},\ucolumnvector{j})\bigr).\]
Moreover, $\mycomponentvector{k}$ and $\mycomponentvector{\ell}$ are linearly independent.
\end{lemma}
\begin{bew}
Let $\boldsymbol{v}\in \R^m$. We calculate the determinant of the $(n+1)\times (n+1)$-matrix
\begin{align}\label{align_matrix_cv}
\begin{pmatrix}
C^\trans A \boldsymbol{v}&C^\trans A U
\end{pmatrix}=
\begin{pmatrix}
C^\trans A \boldsymbol{v}&C^\trans A \ucolumnvector{1}&\ldots &C^\trans A \ucolumnvector{n}
\end{pmatrix}
\end{align}
by expanding along the first column and obtain
\begin{align*} 
\det \begin{pmatrix}
C^\trans A \boldsymbol{v}&C^\trans A U
\end{pmatrix}=\sum_{i=0}^n\bigl(C^\trans A \boldsymbol{v}\bigr)_{i+1}\mysamesign{i}=\sum_{i=0}^n B(\mycomponentvector{i},\boldsymbol{v})\ko\mysamesign{i}.
\end{align*}
For $\boldsymbol{v}=\ucolumnvector{j}$ the determinant vanishes, as (\ref{align_matrix_cv}) has two identical columns, i.\,e.
\begin{align*}
\sum_{i=0}^n B(\mycomponentvector{i},\ucolumnvector{j})\ko\mysamesign{i}=0.
\end{align*}
For any $U\in \mathcal{C}_A$, we cannot have $B(\mycomponentvector{i},\ucolumnvector{j})=0$ for all $i\in \lbrace 0,\ldots,n\rbrace$, since otherwise the $j$-th row of $U^\trans A \widetilde{C}_k$
is zero, which implies $\mysamesign{k}=0$ for all $k\in \lbrace 0,\ldots,n\rbrace$. 
If $B(\mycomponentvector{i},\ucolumnvector{j})>0$ for all $i\in \lbrace 0,\ldots,n\rbrace$, the expression $\sum_{i=0}^n B(\mycomponentvector{i},\ucolumnvector{j})\mysamesign{i}$ would be strictly positive (resp. negative) since $\mysamesign{k}\geq 0$ (resp. $\mysamesign{k}\leq 0$) for all $k\in \lbrace 0,\ldots,n\rbrace$, where at least one inequality is strict.
With the same argument, we exclude the case $B(\mycomponentvector{i},\ucolumnvector{j})<0$ for all $i\in \lbrace 0,\ldots,n\rbrace$.
Hence, there exist $k,\ell\in\lbrace 0,\ldots,n\rbrace$ with $\sgn\bigl(B(\mycomponentvector{k},\ucolumnvector{j})\bigr)\neq \sgn\bigl(B(\mycomponentvector{\ell},\ucolumnvector{j})\bigr)$.

We note that any linearly dependent column vectors $\mycomponentvector{k}$ and $\mycomponentvector{\ell}$ admit the same sign: Let $\lambda \in \R$ such that $\mycomponentvector{\ell}=\lambda \mycomponentvector{k}$. Then $B(\mycomponentvector{k},\mycomponentvector{\ell})=\lambda B(\mycomponentvector{k},\mycomponentvector{k})$ holds. As we have $B(\mycomponentvector{k},\mycomponentvector{\ell})<0$ and $B(\mycomponentvector{k},\mycomponentvector{k})<0$ for $\mycomponentvector{k},\mycomponentvector{\ell}\in\comp$, the factor $\lambda$ is strictly positive, hence $\sgn\bigl(B(\mycomponentvector{\ell},\ucolumnvector{j})\bigr)=\sgn\bigl(\lambda B(\mycomponentvector{k},\ucolumnvector{j})\bigr)= \sgn\bigl(B(\mycomponentvector{k},\ucolumnvector{j})\bigr).$
\end{bew}
We use Lemma \ref{lemma_pos_def} and Lemma \ref{lemma_sign} to construct holomorphic theta series.
\begin{satz}\label{proposition_holomorphic}
The series defining $\mytheta{f}$
is absolutely convergent. Moreover, $\mytheta{f}$ is a holomorphic function in $Z\in \H_n$.
\end{satz}
\begin{bew}
We have $f(U)=0$ for $U\notin \mathcal{C}_A$ as noted above. For $U\in \mathcal{C}_A$ we obtain:
\begin{align*}
f(U)=\begin{cases}
\begin{alignedat}{2}
\phantom{(-1)^n}1&\quad\text{if } \mysamesign{i}>0\quad\text{for all }0\leq i\leq n,\\[5pt] 
-1&\quad\text{if }\mysamesign{i}<0\quad\text{for all }0\leq i\leq n,\\[5pt] 
\phantom{-}2^{k-n-1}&\quad\text{if }\mysamesign{i}>0\quad\text{for $k$ values in }0\leq i\leq n\text{ and }\mysamesign{i}=0\text{ otherwise,}\\[5pt] 
-2^{k-n-1}&\quad\text{if }\mysamesign{i}<0\quad\text{for $k$ values in }0\leq i\leq n\text{ and }\mysamesign{i}=0\text{ otherwise.}
\end{alignedat}
\end{cases}
\end{align*}
Since only the values $U\in \mathcal{C}_A$ contribute non-vanishing terms, we split up $\mathcal{C}_A$ in smaller components. In each component, the expression $\tr (U^\trans A U)$ is bounded from below by a positive definite quadratic form built from the quadratic forms $Q^+_{k,\ell}$ that were introduced in Lemma \ref{lemma_pos_def}. For a fixed $\ucolumnvector{j}$, we apply Lemma \ref{lemma_sign} and take $k,\ell$ with $\sgn\bigl(B(\mycomponentvector{k},\ucolumnvector{j})\bigr)\neq \sgn\bigl(B(\mycomponentvector{\ell},\ucolumnvector{j})\bigr)$, so $B(\mycomponentvector{k},\ucolumnvector{j})\ko B(\mycomponentvector{\ell},\ucolumnvector{j})\leq 0$. Since $\mycomponentvector{k},\mycomponentvector{\ell}\in\comp$ and these vectors are linearly independent, we have 
\[\frac{B(\mycomponentvector{k},\mycomponentvector{\ell})}{4Q(\mycomponentvector{k})Q(\mycomponentvector{\ell})-B(\mycomponentvector{k},\mycomponentvector{\ell})^2}>0.\]
Thus, $Q(\ucolumnvector{j})\geq \Qkl (\ucolumnvector{j})$.
Considering $\Quadratic (U)=\sum_{j=1}^nQ(\ucolumnvector{j})$, the quadratic forms $\Qkl$ give a lower bound for each $\ucolumnvector{j}$ and thus also a bound for $U\in\mathcal{C}_A$. Note that in general we have to take different forms $\Qkl$ for each column vector.
As $Y\in \R^{n\times n}$ is a symmetric positive definite matrix, the square root $Y^{1/2}\in \R^{n\times n}$ is uniquely determined and positive definite. The set $\mathcal{C}_A$ is invariant under the substitution $U\mapsto \breve{U}=UY^{1/2}$ and we find for every column $\boldsymbol{\breve{u}_j}$ of $\breve{U}$ a lower bound in terms of a positive definite quadratic form as before.
Hence,
\[\big|\exp \bigl( \pi i \tr (U^\trans AUZ)\bigr)\big|=\exp \bigl(- \pi \tr (U^\trans AUY)\bigr)=\exp \Bigl(-2 \pi \sum\limits_{j=1}^nQ(\boldsymbol{\breve{u}_j})\Bigr)\] 
is bounded from above by the sum of positive definite quadratic forms $\Qkl$ that we choose for all column vectors $\boldsymbol{\breve{u}_j}$ independently.
We split up the sum over $U\in \mathcal{C}_A \cap \Z^{m\times n}$ in finitely many sets, according to the quadratic forms $\Qkl$ that give a lower bound for $\Quadratic(U)$. 
Thus, the series is absolutely convergent.
Since $f$ is locally constant (the points of discontinuity are given by the matrices in $\mathcal{H}+\Z^{m\times n}$ with $\mysamesign{i}=0$ for some $i\in\lbrace 0,\ldots,n\rbrace$), the series $\mytheta{f}$ is holomorphic in $Z$.
\end{bew}
\begin{remark}
We give a specific formula for $f$ here, but we can replace $f$ by any locally constant function that is zero for $U\notin\mathcal{C}_A$ to obtain a holomorphic theta series.
\end{remark}
This shows part (i) of Theorem \ref{maintheorem}.
In the following section, we construct certain functions $g$ (depending on the choice of $C$) such that $\mytheta{g}$ has modular transformation properties. 
We will see that  $g(UY^{1/2})\rightarrow f(U)$ almost everywhere for $Y\rightarrow\infty$.
\section{Siegel theta series with modular transformation behavior}\label{section_modular}
In \cite{roehrig} we have constructed theta series $\mytheta{g}$ that transform like modular forms by considering a certain family of functions $g$, which we described here in Theorem \ref{theorem_solutions}. A crucial attribute of these functions is that we can split up $g$ in two factors where one depends on a subspace of $\R^{m\times n}$, on which the quadratic form is positive semi-definite, and the other on a subspace where the form is negative semi-definite.
In the following, we first determine explicitly how we split up the quadratic form for matrices of signature $(m-1,1)$.
Then we show that we can apply the result of \cite{roehrig} to deduce the modular transformation behavior of the theta series.

In the next lemma, we show that $\Quadratic$ is positive semi-definite on $\Usenk=(\ucolumnvectorsenk{1}\ldots\ucolumnvectorsenk{n})$ and negative semi-definite on $\Uc=(\ucolumnvectorc{1}\ldots\ucolumnvectorc{n})$.
\begin{lemma}\label{lemma_split}
For $\Quadratic :\R^{m\times n}\longrightarrow \R$, we have the decomposition $\Quadratic=\Quadraticplus+\Quadraticminus$, where $\Quadraticplus$ is positive semi-definite and $\Quadraticminus$ negative semi-definite. Moreover, $\Quadraticplus(U)=\Quadratic(\Usenk)$ and $\Quadraticminus(U)=\Quadratic(\Uc)$.
\end{lemma}
\begin{bew}
By the definition of $\Aneg$ in \eqref{align_split}, we immediately obtain
\begin{align}\label{align_Qminus}
\Quadraticminus(U)=\tr  \biggl(\frac{U^\trans A\mycomponentvectorgeneral\mycomponentvectorgeneral^\trans AU}{4Q(\mycomponentvectorgeneral)}\biggr)=\frac1{4Q(\mycomponentvectorgeneral)} \sum\limits_{j=1}^nB(\mycomponentvectorgeneral,\ucolumnvector{j})^2 \leq 0 \quad\text{for all }U\in \R^{m\times n}.
\end{align}
Since $\Apos=A-\Aneg$, we have
\begin{align}\label{align_Qplus}
\begin{split}
\Quadraticplus(U)&=\frac12\tr (U^\trans A U)-\tr \biggl(\frac{U^\trans A\mycomponentvectorgeneral\mycomponentvectorgeneral^\trans AU }{4Q(\mycomponentvectorgeneral)}\biggr)\\
&=\sum\limits_{j=1}^{n} Q(\ucolumnvector{j})-\frac{1}{4Q(\mycomponentvectorgeneral)}\sum\limits_{j=1}^n B(\mycomponentvectorgeneral,\ucolumnvector{j})^2\\
&=\sum\limits_{j=1}^n \frac{4Q(\mycomponentvectorgeneral)Q(\ucolumnvector{j})-B(\mycomponentvectorgeneral,\ucolumnvector{j})^2}{4Q(\mycomponentvectorgeneral)}.
\end{split}
\end{align}
The numerator of each summand represents the determinant of the Gram matrix
\begin{align*}
\left(\begin{matrix}
2Q(\ucolumnvector{j})&B(\mycomponentvectorgeneral,\ucolumnvector{j})\\
B(\mycomponentvectorgeneral,\ucolumnvector{j})&2Q(\mycomponentvectorgeneral)
\end{matrix}\right).
\end{align*}
For linearly independent vectors $\mycomponentvectorgeneral$ and $\ucolumnvector{j}$, the quadratic form $Q$ has signature $(1,1)$ on $\spann \lbrace \mycomponentvectorgeneral ,\ucolumnvector{j}\rbrace$, i.\,e. the Gram matrix has negative determinant. For linearly dependent vectors $\mycomponentvectorgeneral$ and $\ucolumnvector{j}$, we obtain zero.  As $Q(\mycomponentvectorgeneral)<0$, we thus see, using \eqref{align_Qplus}, that $\Quadraticplus(U)\geq 0$ for all $U\in \R^{m\times n}$. Note that $\Quadraticplus(U)=0$ holds if and only if every column of $U$ is a multiple of $\mycomponentvectorgeneral$.

The negative semi-definite part $\Quadraticminus$ only depends on $\Uc$. This follows immediately when we use the identity $Q(\ucolumnvectorc{j})=\frac{B(\mycomponentvectorgeneral,\ucolumnvector{j})^2}{4Q(\mycomponentvectorgeneral)}:$
\begin{align*}
\Quadraticminus(U)&=\frac1{4Q(\mycomponentvectorgeneral)} \sum\limits_{j=1}^nB(\mycomponentvectorgeneral,\ucolumnvector{j})^2
=\sum\limits_{j=1}^{n} Q(\ucolumnvectorc{j})
=\frac12\tr \bigl((\Uc)^\trans A \Uc\bigr)
=\Quadratic(\Uc)
\end{align*} 
Then the positive semi-definite quadratic form $\Quadraticplus$ only depends on $\Usenk$, i.\,e. the part of $U$ that is perpendicular to $\mycomponentvectorgeneral$,
because $Q(\ucolumnvector{j})=Q(\ucolumnvectorsenk{j})+Q(\ucolumnvectorc{j})$. Thus, we have  $\Quadratic=\Quadraticplus+\Quadraticminus$ with $\Quadraticplus(U)=\Quadratic(\Usenk)$ and $\Quadraticminus(U)=\Quadratic(\Uc)$.
\end{bew}
Before we prove the remaining parts (ii) and (iii) of Theorem \ref{maintheorem}, we recall Zwegers' construction \cite{zwegers} for the case $n=1$, as that makes clear how we choose the set-up for higher genus $n$. We merely give the function $h$ that defines $\mytheta{h}$ here:
\begin{defi}[{\cite[Definition 2.1]{zwegers}}]\label{defi_zwegers}
Let $\mycomponentvector{0},\mycomponentvector{1}\in\comp \subset \R^m$ and define
\begin{align*}
h(\uvector)=h^{\mycomponentvector{0},\mycomponentvector{1}}(\uvector):=E\Bigl(\frac{B(\mycomponentvector{0},\boldsymbol{u})}{\sqrt{-Q(\mycomponentvector{0})}}\Bigr)-E\Bigl(\frac{B(\mycomponentvector{1},\boldsymbol{u})}{\sqrt{-Q(\mycomponentvector{1})}}\Bigr).
\end{align*}
\end{defi}
Kudla \cite{kudla13} and Livinsky \cite{livinsky} showed that the corresponding theta series $\mytheta{h}$ can be constructed as integrals of the theta forms introduced by Kudla and Millson \cite{KM86,KM87}. We do something similar here, that is we show in the next lemma that $h$ is obtained by integrating a certain 1-form over $S_1$. To this end, we have to choose an explicit parameterization of $S_1$, here we consider
\[S_1=\lbrace t\ko\mycomponentvector{0}+(1-t)\ko\mycomponentvector{1}\mid t\in [0,1]\rbrace.\]
\begin{lemma}\label{lemma_differentiate}
We can write $h$ as
\[h(\uvector)=2\int_{S_1} \exp\Bigl(\pi \frac{B(\mycomponentvectorgeneral,\uvector)^2}{Q(\mycomponentvectorgeneral)}\Bigr)\ko \frac{B(\usenk,d\mycomponentvectorgeneral)}{\sqrt{-Q(\mycomponentvectorgeneral)}}\quad\text{with } d\mycomponentvectorgeneral=(dc_1,\ldots,dc_m)^\trans .\]
\end{lemma}
\begin{proof}
Note that
\[\frac{\partial}{\partial \mycomponentvectorgeneral} B(\mycomponentvectorgeneral,\boldsymbol{u})=A\boldsymbol{u}\quad\text{and}\quad \frac{\partial}{\partial \mycomponentvectorgeneral} Q(\mycomponentvectorgeneral)=A\mycomponentvectorgeneral.\]
Since $\uvector=\usenk +\frac{B(\mycomponentvectorgeneral,\uvector)}{2Q(\mycomponentvectorgeneral)}\mycomponentvectorgeneral$, we obtain
\begin{align*}
\frac{\partial}{\partial \mycomponentvectorgeneral} \Bigl(\frac{B(\mycomponentvectorgeneral,\uvector)}{\sqrt{-Q(\mycomponentvectorgeneral)}}\Bigr)
=\frac{A\uvector \sqrt{-Q(\mycomponentvectorgeneral)}+A\mycomponentvectorgeneral\ko B(\mycomponentvectorgeneral,\uvector) \bigl(2\sqrt{-Q(\mycomponentvectorgeneral)}\bigr)^{-1}}{-Q(\mycomponentvectorgeneral)}
=\frac{A\usenk}{\sqrt{-Q(\mycomponentvectorgeneral)}}.
\end{align*}
Hence,
\[\frac{\partial}{\partial \mycomponentvectorgeneral} E\Bigl(\frac{B(\mycomponentvectorgeneral,\boldsymbol{u})}{\sqrt{-Q(\mycomponentvectorgeneral)}}\Bigr)=\frac{A\usenk}{\sqrt{-Q(\mycomponentvectorgeneral)}}\ko E'\Bigl(\frac{B(\mycomponentvectorgeneral,\boldsymbol{u})}{\sqrt{-Q(\mycomponentvectorgeneral)}}\Bigr)=2\ko\frac{A\usenk}{\sqrt{-Q(\mycomponentvectorgeneral)}}\ko\exp\Bigl(\pi \frac{B(\mycomponentvectorgeneral,\boldsymbol{u})^2}{Q(\mycomponentvectorgeneral)}\Bigr),\]
so the total differential of $E$ with regard to $\mycomponentvectorgeneral$ is the exact 1-form
\begin{align}\label{align_differentialE}
dE\Bigl(\frac{B(\mycomponentvectorgeneral,\uvector)}{\sqrt{-Q(\mycomponentvectorgeneral)}}\Bigr)=2\ko\exp\Bigl(\pi \frac{B(\mycomponentvectorgeneral,\uvector)^2}{Q(\mycomponentvectorgeneral)}\Bigr)\ko\frac{B(\usenk,d\mycomponentvectorgeneral)}{\sqrt{-Q(\mycomponentvectorgeneral)}}.
\end{align}
We integrate both sides of (\ref{align_differentialE}) over $S_1$ to finish the proof.
\end{proof}
We transfer this construction to Siegel theta series of generic genus $n\in\N$ by considering the integrand of the function $h$ from Lemma \ref{lemma_differentiate} for each column vector $\ucolumnvector{j}$ of $U$ and taking the wedge product over all $j=1,\ldots ,n$ to obtain an (exact) $n$-form that is integrated over the $n$-simplex $S_n$.

Note that this is an explicit realization of the theta forms $\theta_{KM}$ for arbitrary signature valued in closed differential forms that were constructed by Kudla and Millson \cite{KM86,KM87,KM90}. The next proposition is a result that was also shown by Livinsky \cite{livinsky} in his Ph.\,D. thesis (based on an unpublished manuscript by Kudla \cite{kudla13}): he defines $\Theta_{KM}^{\Delta}$ as the integral of the closed $n$-form $\theta_{KM}$ over the simplex $\Delta$ (which is $S_n$ in our notation) and thus constructs a non-holomorphic Siegel modular form. 

We make a similar construction but instead of using the connection to the theta forms $\theta_{KM}$, we show that we obtain functions that we already know from \cite{roehrig}, which also shows that the Siegel theta series that we obtain are modular.
\begin{satz}\label{proposition_modular}
The theta series $\mytheta{g}$ transforms like a Siegel modular form of weight $m/2$.
\end{satz}
\begin{bew}
We recall that the integrand of $g$ is
\begin{align*}
\bigl(-Q(\mycomponentvectorgeneral)\bigr)^{-n/2}\ko\exp\bigl(2\pi \tr (U^\trans \Aneg U)\bigr)\ko\bigwedge_{j=1}^n  B(\ucolumnvectorsenk{j},d\mycomponentvectorgeneral).
\end{align*}
Now we first consider the wedge product
and write the bilinear forms as sums. Using the distributivity of the wedge product, we obtain
\begin{align*}
\bigwedge_{j=1}^n\biggr(\sum\limits_{k_j=1}^m (A\ucolumnvectorsenk{j})_{k_j} dc_{k_j}\biggr)
&=\sum\limits_{k_1,\ldots, k_n=1}^m (A\ucolumnvectorsenk{1})_{k_1}\cdots (A\ucolumnvectorsenk{n})_{k_n} dc_{k_1}\wedge \ldots \wedge dc_{k_n}.
\end{align*}
As we have $dc_{k_1}\wedge \ldots \wedge dc_{k_n}=0$ if $k_i=k_j$ for any $i\neq j$ and  $dc_{k_1}\wedge \ldots \wedge dc_{k_n}=\sgn (\sigma) dc_{k_{\sigma(1)}}\wedge \ldots \wedge dc_{k_{\sigma(n)}}$ for any permutation in the symmetric group $\sigma \in S_n$, this expression equals
\begin{align*}
\sum\limits_{1\leq k_1<k_2<\ldots< k_n\leq m}\ko\biggl(\sum\limits_{\sigma\in S_n} \sgn (\sigma) \prod\limits_{j=1}^{n}(A\ucolumnvectorsenk{j})_{k_{\sigma(j)}}\biggr)\ko dc_{k_1}\wedge \ldots \wedge dc_{k_n}.
\end{align*}
We observe that the part in brackets is the determinant of the matrix 
\begin{align}\label{align_submatrix}
\bigl((A\ucolumnvectorsenk{j})_{k_i}\bigr)_{ij}=
\begin{pmatrix}
(A\ucolumnvectorsenk{1})_{k_1}&(A\ucolumnvectorsenk{2})_{k_1}&\cdots & (A\ucolumnvectorsenk{n})_{k_1}\vspace{4pt}\\
(A\ucolumnvectorsenk{1})_{k_2}&(A\ucolumnvectorsenk{2})_{k_2}&\cdots & (A\ucolumnvectorsenk{n})_{k_2}\\
\vdots&\vdots& &\vdots\\
(A\ucolumnvectorsenk{1})_{k_n}&(A\ucolumnvectorsenk{2})_{k_n}&\cdots & (A\ucolumnvectorsenk{n})_{k_n}
\end{pmatrix}\quad (1\leq i \leq n,\ko1\leq j \leq n),
\end{align}
which is a square submatrix of maximal size of $A\Usenk\in \R^{m\times n}$. We use multi-index notation to rewrite this. Let 
\[K:=\lbrace\kvector =(k_1,\ldots,k_n)\in \N^n\mid 1\leq k_1<k_2<\ldots< k_n\leq m\rbrace\]
and for $\mathcal{A}\in \R^{m\times n}$ let us denote by $\mathcal{A}_{\kvector}\in \R^{n\times n}$ the square submatrix that consists of the rows determined by $\kvector$. So \eqref{align_submatrix} can be written as $(A\Usenk)_{\kvector}$ for $\kvector \in K$.
Obviously, $P_{\kvector}(U):=\det\bigl((A\Usenk)_{\kvector}\bigr)$ has the homogeneity property
$P_{\kvector}(UN)=\det N\cdot P_{\kvector}(U)$ for $N\in \C^{n\times n}$, with the degree of homogeneity being 1. Applying the Laplacian $\laplacian_{\Amaj}$ means differentiating twice with regard to every row, so a function with the determinant-like structure of $P_{\kvector}$ and degree 1 will necessarily vanish under this operator, i.\,e. $(\laplacian_{\Amaj}) P_{\kvector}=0$. So $P_{\kvector}$ is harmonic and we simply have $\exp\bigl( -\laplacian_{\Amaj}/8\pi\bigr)\bigl(P_{\kvector}(U)\bigr)=P_{\kvector}(U)$. 
Further, we note that $P_{\kvector}$ only depends on a subspace of $\R^{m\times n}$ where the quadratic form $\Quadratic$ is positive semi-definite.

By Lemma \ref{lemma_split}, we have $$\exp\bigl(2\pi \tr (U^\trans A^{-}U)\bigr)=\exp\bigl(4\pi \Quadraticminus(U)\bigr)=\exp\bigl(4\pi \ko \Quadratic (\Uc )\bigr),$$ so the exponential factor solely depends on a subspace of $\R^{m\times n}$ on which $\Quadratic$ is negative semi-definite.
Thus, the integrand of $g$ has the form
\begin{align}\label{align_integrand}
\bigl(-Q(\mycomponentvectorgeneral)\bigr)^{-n/2}\ko\exp\bigl(4\pi \ko \Quadratic (\Uc )\bigr)\ko\sum_{\kvector \in K}  \det\bigl((A\Usenk)_{\kvector}\bigr)\ko d\mycomponentvectorgeneral_{\kvector}\quad\text{with }d\mycomponentvectorgeneral_{\kvector}=dc_{k_1}\wedge \ldots \wedge dc_{k_n},
\end{align}
which is a function as described in Theorem \ref{theorem_solutions}. Using the notation of this theorem, we have a function where $P_r$ is a harmonic polynomial of degree $\alpha=1$ and $P_s\equiv 1$ (and so has degree $\beta=0$). So the associated theta series transforms like a Siegel modular form of weight $m/2$.

In Remark \ref{rem_Sn}, we observed that $S_n\subset \comp$ holds and that $S_n$ is compact in $\R^m$. Hence, the points in $S_n$ do not accumulate near the boundary of $\comp$, i.\,e. where $Q(\mycomponentvectorgeneral)$ is almost zero. So the summands of the theta series associated with \eqref{align_integrand} are rapidly decaying functions and we can integrate termwise over $S_n$ to obtain $\mytheta{g}$. The modular transformation properties are preserved as they are independent of the choice of $C$. Thus $\mytheta{g}$ is well-defined and transforms as a modular Siegel theta series of weight $m/2$.
\end{bew}
In the definition of $\mytheta{g}$ we consider $g(UY^{1/2})$ instead of $g(U)$, where $Y$ denotes the imaginary part of $Z$. We are interested in the behavior of the theta series for large values of $Y$. For $n=1$, it is clear that we consider the imaginary part $y\in \R_{>0}$ as large, when $y\rightarrow \infty$. We recall that in Definition \ref{defi_zwegers} the error function $E$ was considered, where $E(xy^{1/2})\rightarrow \sgn (x)$ for $y\rightarrow \infty$.
For arbitrary genus $n\in \N$, we gave a generalizing definition of $Y\rightarrow\infty$ in Definition \ref{defi_infinity}. We show in the next proposition that $g(UY^{1/2})$ asymptotes to the locally constant function $f(U)$ for $Y\rightarrow\infty$.
\begin{satz}\label{proposition_completion}
We have $g(UY^{1/2})\rightarrow f(U)$ almost everywhere for $Y\rightarrow\infty$.
\end{satz}
\begin{bew}
Using identity \eqref{align_Qminus}, we can write the integrand of $g$ as
\[\exp\Bigl(\frac{\pi}{Q(\mycomponentvectorgeneral)}\ko \sum\limits_{j=1}^{n}B(\mycomponentvectorgeneral,\ucolumnvector{j})^2\Bigr)\ko\bigwedge_{j=1}^n  \frac{B(\ucolumnvectorsenk{j},d\mycomponentvectorgeneral)}{\sqrt{-Q(\mycomponentvectorgeneral)}}.\]
We substitute
\begin{align*}
\vvector:=\frac{U^\trans A\mycomponentvectorgeneral}{\sqrt{-Q(\mycomponentvectorgeneral)}}=\frac{1}{\sqrt{-Q(\mycomponentvectorgeneral)}}\left(\begin{array}{c}
B(\mycomponentvectorgeneral,\ucolumnvector{1})\\
\vdots\\
B(\mycomponentvectorgeneral,\ucolumnvector{n})
\end{array}\right).
\end{align*}
In the proof of Lemma \ref{lemma_differentiate} we have shown  that \[\frac{\partial}{\partial \mycomponentvectorgeneral} \Bigl(\frac{B(\mycomponentvectorgeneral,\uvector)}{\sqrt{-Q(\mycomponentvectorgeneral)}}\Bigr)=\frac{A\usenk}{\sqrt{-Q(\mycomponentvectorgeneral)}},\] so the differential of $v_j= \frac{B(\mycomponentvectorgeneral,\ucolumnvector{j})}{\sqrt{-Q(\mycomponentvectorgeneral)}}$ with respect to $\mycomponentvectorgeneral$ is easily seen to be
$dv_j=\frac{B(\ucolumnvectorsenk{j},d\mycomponentvectorgeneral)}{\sqrt{-Q(\mycomponentvectorgeneral)}}$.
We can thus write $g$ as
\begin{align*}
\int\limits_{X_n} \exp \bigl(-\pi (v_1^2+\ldots +v_n^2)\bigr)\ko dv_1\wedge \ldots\wedge dv_n,
\end{align*}
where
\begin{align*}
X_n:=\left\lbrace \frac{t_0\mysamesigncolumnvector{0}+t_1\mysamesigncolumnvector{1}+\ldots +t_n\mysamesigncolumnvector{n}}{\sqrt{-Q(t_0\mycomponentvector{0}+t_1\mycomponentvector{1}+\ldots +t_n\mycomponentvector{n})}}\hspace{3pt}\bigg\vert\hspace{3pt} t_i\geq 0\text{ for all }0\leq i \leq n\quad\text{and}\quad\sum\limits_{i=0}^n t_i=1 \right\rbrace\quad\bigl(\mysamesigncolumnvector{i}=U^\trans A \mycomponentvector{i}\bigr).
\end{align*}
We describe this integral depending on whether $U$ is in $\mathcal{C}_A$ or not. If
\begin{align}\label{align_determinant_x_i}
\det\begin{pmatrix}
1&1&\ldots &1\\
\mysamesigncolumnvector{0}&\mysamesigncolumnvector{1}&\ldots &\mysamesigncolumnvector{n}
\end{pmatrix}=\sum_{i=0}^n \mysamesign{i}    
\end{align}
vanishes, we have $U\notin \mathcal{C}_A$, as either all $\mysamesign{i}$ vanish or there occur $\mysamesign{i}$ and $\mysamesign{j}$ that have different signs. Moreover, considering the left-hand side of \eqref{align_determinant_x_i}, the vectors $\mysamesigncolumnvector{0},\ldots,\mysamesigncolumnvector{n}$ are linearly dependent, and so the vertices $\mysamesigncolumnvector{i}/\sqrt{-Q(\mycomponentvector{i})}\,(0\leq i \leq n)$ form a simplex whose dimension is strictly lower than $n$. But when integrating an $n$-form over this simplex, the integral takes the value zero.

If \eqref{align_determinant_x_i} does not vanish, we have $\boldsymbol{0}\in X_n$ if and only if $U\in \mathcal{C}_A$, which is a direct consequence of Cramer's rule: the system of $n+1$ linear equations
\begin{align*}
\begin{pmatrix}
1&1&\ldots &1\\
\mysamesigncolumnvector{0}&\mysamesigncolumnvector{1}&\ldots &\mysamesigncolumnvector{n}
\end{pmatrix}
\begin{pmatrix}
t_0\\
t_1\\
\vdots\\
t_n
\end{pmatrix}=
\begin{pmatrix}
1\\
0\\
\vdots\\
0
\end{pmatrix}
\end{align*}
has a unique solution $(t_0,t_1,\ldots,t_n)\in \R^{n+1}$ with $t_i\geq 0$ for all $0\leq i\leq n$ if and only if \[\mysamesign{i}=(-1)^{i}\det \bigl(\mysamesigncolumnvector{0}\ldots\widehat{\mysamesigncolumnvector{i}}\ldots \mysamesigncolumnvector{n})\quad (0\leq i\leq n)\] are all non-negative or all non-positive. If $\mysamesign{i}=0$ for any $i\in\lbrace 0,\ldots,n\rbrace$, the zero vector is located on the boundary of $X_n$.

For $U\in \mathcal{C}_A$, we also determine the orientation of $X_n$ in $\R^n$. Let the canonical basis $\{\boldsymbol{e_1},\ldots,\boldsymbol{e_n}\}$ of $\R^n$ represent the equivalence class of positive orientations in $\R^n$. We consider the unit simplex with the vertices $\boldsymbol{0},\boldsymbol{e_1},\ldots,\boldsymbol{e_n}$ that we can also write as
\[T_n:=\bigg\lbrace (t_1,\ldots,t_n)\in\R^{n}\hspace{3pt}\Big\vert\hspace{3pt} t_i\geq 0\text{ for all }1\leq i \leq n\quad\text{and}\quad \sum_{i=1}^n t_i\leq 1\bigg\rbrace.\]
This simplex has a positive orientation in $\R^{n}$. When we consider a diffeomorphism between $T_n$ and any other $n$-simplex in $\R^n$, we can thus determine whether this orientation is preserved (then this simplex also carries the positive orientation in $\R^n$) or reversed. Setting
\[X_n'=\left\lbrace \mysamesigncolumnvector{0}+\sum_{i=1}^n t_i (\mysamesigncolumnvector{i}-\mysamesigncolumnvector{0})\hspace{3pt}\bigg\vert\hspace{3pt} t_i\geq 0\text{ for all }1\leq i \leq n\quad\text{and}\quad\sum\limits_{i=1}^n t_i\leq 1 \right\rbrace,\]
we consider the diffeomorphism \[\varphi:T_n\longrightarrow X_n',\quad \tvector=(t_1,\ldots,t_n)\mapsto \mysamesigncolumnvector{0}+\sum_{i=1}^n t_i (\mysamesigncolumnvector{i}-\mysamesigncolumnvector{0}).\]
We denote the differential of $\varphi$ by $D\varphi$ and observe that $\det \bigl(D\varphi(\tvector)\bigr)$ is independent of $\tvector$ and -- adding an extra column and row -- equals
\[
\det \begin{pmatrix}\mysamesigncolumnvector{1}-\mysamesigncolumnvector{0}&\mysamesigncolumnvector{2}-\mysamesigncolumnvector{0}&\ldots& \mysamesigncolumnvector{n}-\mysamesigncolumnvector{0}\end{pmatrix}=\det \begin{pmatrix}1&0&\ldots &0\\
\mysamesigncolumnvector{0}&\mysamesigncolumnvector{1}-\mysamesigncolumnvector{0}&\ldots& \mysamesigncolumnvector{n}-\mysamesigncolumnvector{0}\end{pmatrix}.
\]
When adding the first column to each of the other columns, we observe that we obtain the expression in \eqref{align_determinant_x_i}, so the determinant of $D\varphi$ is strictly positive if $\mysamesign{i}\geq 0$ for all $i\in \{0,\ldots,n\}$ and strictly negative if $\mysamesign{i}\leq 0$ (in both cases at least one of the inequalities is strict, as $U\in \mathcal{C}_A$).
Up to a normalization, this is $X_n$, which we can write as
\[X_n=\left\lbrace \frac{\mysamesigncolumnvector{0}+\sum_{i=1}^n t_i (\mysamesigncolumnvector{i}-\mysamesigncolumnvector{0})}{\sqrt{-Q\bigl(\mycomponentvector{0}+\sum_{i=1}^n t_i (\mycomponentvector{i}-\mycomponentvector{0})\bigr)}}\hspace{3pt}\bigg\vert\hspace{3pt} t_i\geq 0\text{ for all }1\leq i \leq n\quad\text{and}\quad\sum\limits_{i=1}^n t_i\leq 1 \right\rbrace.\]
Since $S_n\subset \comp$ as observed in Remark \ref{rem_Sn}, we have $\sqrt{-Q\bigl(\mycomponentvector{0}+\sum_{i=1}^n t_i (\mycomponentvector{i}-\mycomponentvector{0})\bigr)}>0$, so the orientation of $X_n$ and $X_n'$ agree.
We make a final substitution, which does not change the orientation: we substitute $U\mapsto UY^{1/2}$ for $Y=\Im Z$, which means that we set $\boldsymbol{\breve{\xvector}_i}=Y^{1/2}\mysamesigncolumnvector{i}$ and integrate over
\[\breve{X}_n=\left\lbrace \frac{t_0\boldsymbol{\breve{\xvector}_0}+t_1\boldsymbol{\breve{\xvector}_1}+\ldots +t_n\boldsymbol{\breve{\xvector}_n}}{\sqrt{-Q(t_0\mycomponentvector{0}+t_1\mycomponentvector{1}+\ldots +t_n\mycomponentvector{n})}}\hspace{3pt}\bigg\vert\hspace{3pt} t_i\geq 0\text{ for all }0\leq i \leq n\quad\text{and}\quad\sum\limits_{i=0}^n t_i=1 \right\rbrace\]
instead of $X_n$.
As $Y$ is positive definite and symmetric, so is its square root $Y^{1/2}$. Since $\det Y^{1/2}> 0$, the property that $\boldsymbol{0}\in \breve{X}_n$ if and only if $U\in \mathcal{C}_A$ is still maintained. While the integrand is still the same as before, the set $\breve{X}_n$ also depends on $Y^{1/2}$. Now we determine how  $\breve{X}_n$ changes when we consider $Y\rightarrow \infty$, so first we observe which property the square root of $Y$ necessarily satisfies for large $Y$.

We have defined $Y\rightarrow \infty$ as $Y_{jj}\rightarrow \infty$ for all $j\in\lbrace 1,\ldots ,n\rbrace$. By definition of $Y^{1/2}$ we have $Y=Y^{1/2}\cdot Y^{1/2}$, which means that the $j$-th diagonal entry is $Y_{jj}=\sum_{\nu=1}^{n}\bigl((Y^{1/2})_{j\nu}\bigr)^2$. So when $Y\rightarrow \infty$, in each row of $Y^{1/2}$ the absolute value of at least one entry tends to infinity. 
Thus, every entry of $\boldsymbol{\breve{\xvector}_i}$ tends to $\pm \infty$ and as the vertices of $\breve{X}_n$ expand outwards, the object we obtain depends on the location of $\boldsymbol{0}$ in relation to this simplex.

For any $U$ such that $\boldsymbol{0}\notin \breve{X}_n$, the set is shifted away from $\boldsymbol{0}$. As the integrand decays fast for large values of $\vvector\in \R^n$, the value of the integral and thus the whole expression $g(UY^{1/2})$ tends to zero. By the definition of $f$, we also have $f(U)=0$ for $U\notin \mathcal{C}_A$, so $g(UY^{1/2})\rightarrow f(U)$ for $Y\rightarrow \infty$ here.

If $\boldsymbol{0}$ is an interior point of $\breve{X}_n$, the simplex asymptotes to $\R^n$. But we know that we have 
\[\int\limits_{\R^n} \exp \bigl(-\pi (v_1^2+\ldots +v_n^2)\bigr)\ko dv_1\wedge \ldots\wedge dv_n= 1,\]
when the canonical basis $\{\boldsymbol{e_1},\ldots,\boldsymbol{e_n}\}$ fixes a positive orientation in $\R^n$.
Thus $g(UY^{1/2})\rightarrow \pm 1$ for $Y\rightarrow\infty$, where the sign depends on the orientation of the original $n$-simplex $X_n$ in $\R^n$ that we described above: the sign is positive if $\mysamesign{i}\geq 0$ for all $1\leq i\leq n$ and negative if $\mysamesign{i}\leq 0$ for all $1\leq i\leq n$.
When $\mysamesign{i}$ are strictly positive or negative for all $i\in\{0,\ldots,n\}$, this is exactly the definition of $f$ for $U\in \mathcal{C}_A$.
So $g(UY^{1/2})\rightarrow f(U)$ for $Y\rightarrow\infty$ almost everywhere. The limit of $g(UY^{1/2})$ may differ from $f(U)$ when $\boldsymbol{0}$ is a boundary point, but the values of $U$ for which this holds form a null set in $\R^{m\times n}$.
\end{bew}
In the following remark, we give a short description of $g(UY^{1/2})$ for $Y\rightarrow\infty$ for the case that $\boldsymbol{0}$ is a boundary point of $\breve{X}_n$.
\begin{remark}\label{rem_solidangle}
Let $1\leq n'\leq n$. When for exactly $n'$ values $\mysamesign{i}=0$ holds, the zero vector is located on an $(n-n')$-face of $\breve{X}_n$. If $n'=1$, this is a facet of $\breve{X}_n$, so the value of $g$ approaches $\pm 1/2$. In this case, we actually see that this agrees with the value of $f(U)$, as $\mysamesign{i}=0$ holds for exactly one $i\in\{0,\ldots,n\}$.

For $n'\geq 2$, the area that we obtain by intersecting $\breve{X}_n$ with the $(n'-1)$-dimensional unit-sphere is called the solid angle $\Omega$. Then $g$ asymptotes to $\pm\Omega/A_{n'}$, where $A_{n'}$ is the surface area of the unit-sphere. Again, the sign reflects the orientation of $X_n$ (resp. $\breve{X}_n$) in $\R^n$, as described in the previous proof.

Considering the solid angle that depends on the exact position of $\boldsymbol{0}$ in the simplex, one could determine an exact formula for the holomorphic part of $\mytheta{g}$. However, the resulting function will look extremely complicated (for $n=2$ one could for example use the result in \cite{livinsky}), so we used the holomorphic function $f$, which has a much simpler form.
\end{remark}
In this section, we have thus shown the second and third part of the main theorem: By Proposition \ref{proposition_modular}, the theta series in Theorem \ref{maintheorem}(ii) transforms like a Siegel modular form of genus $n$ and weight $m/2$. Part (iii), giving us the connection between the holomorphic version and the modular version of the theta series, follows by Proposition \ref{proposition_completion}.
\begin{acknowledgments}
The author would like to kindly thank Sander Zwegers for suggesting the topic, as well as providing a lot of insightful advice.
\end{acknowledgments}
\bibliography{References}

\begin{thebibliography}{99}
\bibitem{ABMP18}
S.~Alexandrov, S.~Banerjee, J.~Manschot, and B.~Pioline, \textit{Indefinite theta series and generalized error functions}, Sel. Math. New Ser. \textbf{24} (2018), no.~5, 3927--3972.

\bibitem{borcherds}
R.~Borcherds, \textit{Automorphic forms with singularities on Grassmannians}, Invent. Math. \textbf{132} (1998), no.~3, 491--562.

\bibitem{Fre83}
E.~Freitag, \textit{Siegelsche Modulfunktionen}, Grundlehren der mathematischen Wissenschaften \textbf{254} (1983), Springer-Verlag Berlin.

\bibitem{FK17}
J.~Funke and S.~Kudla, \textit{Mock modular forms and geometric theta functions for indefinite
quadratic forms}, J.\ Phys.\ A:\ Math.\ Theor.\ \textbf{50} (2017), no.~40.

\bibitem{FK19}
\bysame, \textit{On some incomplete theta integrals}, Compos.\ Math.\ \textbf{155} (2019), no.~9, 1711--1746.

\bibitem{GZ98}
L.~G\"ottsche and D.~Zagier, \textit{Jacobi forms and the structure of Donaldson invariants for 4-manifolds with $b_+=1$}, Selecta Math.\ (N.S.) \textbf{4} (1998), no.~1, 69--115.

\bibitem{kudla81}
S.~Kudla, \textit{Holomorphic Siegel modular forms associated to $\operatorname{SO}(n,1)$}, Math.~Ann. \textbf{256} (1981), 517--534.

\bibitem{kudla13}
\bysame, \textit{A note on Zwegers' theta functions}, preprint (2013).

\bibitem{kudla18}
\bysame, \textit{Theta integrals and generalized error functions}, Manuscripta Math. \textbf{155} (2018), no.~3--4, 303--333.

\bibitem{KM86}
S.~Kudla and J.~Millson, \textit{The theta correspondence and harmonic forms. I}, Math.~Ann. \textbf{274} (1986), no.~3, 353--378.

\bibitem{KM87}
\bysame, \textit{The theta correspondence and harmonic forms. II}, Math.~Ann. \textbf{277} (1987), no.~2, 267--314.

\bibitem{KM90}
\bysame, \textit{Intersection numbers of cycles on locally symmetric spaces and Fourier coefficients of holomorphic modular forms in several complex variables}, Publ.~Math.~IHES \textbf{71} (1990), 121--172.

\bibitem{livinsky}
I.~Livinsky, \textit{On the integrals of the Kudla-Millson theta series}, Ph.D. Dissertation (2016), Toronto.

\bibitem{Naz18}
C.~Nazaroglu, \textit{$r$-tuple error functions and indefinite theta series of higher-depth}, Commun. Number Theory Phys. \textbf{12} (2018), no.~3, 581--608.

\bibitem{roehrig}
C.~Roehrig, \textit{Siegel theta series for indefinite quadratic forms}, preprint (2020), arXiv:2009.08230v2.

\bibitem{schoeneberg}
B.~Schoeneberg, \textit{Das Verhalten von mehrfachen Thetareihen bei Modulsubstitutionen}, Math.~Ann. \textbf{116} (1939), no.~1, 511--523.

\bibitem{shimura}
G.~Shimura, \textit{On modular forms of half integral weight}, Ann.~Math. \textbf{97} (1973), no.~3, 440--481.

\bibitem{raum}
M.~Westerholt-Raum, \textit{Indefinite theta series on cones}, preprint (2016), arXiv:1608.08874.

\bibitem{zwegers}
S.~Zwegers, \textit{Mock Theta Functions}, Ph.D. Dissertation (2002), Utrecht.

\end{thebibliography}

\end{document}